\documentclass[11pt]{article}
\usepackage{amssymb}
\usepackage{amsmath}
\usepackage{enumerate}
\usepackage{amscd}
\usepackage{graphics}
\usepackage{latexsym}

\usepackage{epsfig}
\usepackage[latin1]{inputenc}

\usepackage{a4wide}

\title{
\huge Relative parabolicity of zero mean curvature surfaces in $\r^3$ and $\l^3$}
\author{
\Large Isabel Fern\'{a}ndez   
\thanks{Research partially supported by
MCYT-FEDER grant number MTM2004-00160.\newline
}
\and   \Large Francisco J. L\'{o}pez $ ^{\ast}$  }

\newcommand{\df}{ \stackrel{\rm def}{=}}

\def\wb{\mathcal W}

\def\C{\mathcal C}

\def\h{\mathbb{H}}
\def\r{\mathbb{R}}
\def\n{\mathbb{N}}

\def\c{\mathbb{C}}

\def\s{\mathbb{S}}
\def\d{\mathbb{D}}
\def\l{\mathbb{R}_1}

\def\cb{\mathcal{C}}

\def\rb{\mathcal{R}}
\def\sb{\mathcal{S}}

\def\mb{\mathcal{M}}

\newenvironment{proof}{\trivlist
\item[\hskip\labelsep{\em Proof}\,:]}{\hfill{$\Box$}\endtrivlist}

\newtheorem{lemma}{Lemma}[section]
\newtheorem{remark}{Remark}[section]
\newtheorem{theorem}{Theorem}[section]
\newtheorem{corollary}{Corollary}[section]

\newtheorem{definition}{Definition}[section]

\begin{document}

\maketitle

\begin{abstract}
If the Lorentzian norm on a maximal surface in the 3-dimensional Lorentz-Minkowski space $\l^3$ is positive and proper, then the surface is relative parabolic. As a consequence, entire maximal graphs with a closed set of isolated singularities are relative parabolic.

Furthermore, maximal  and minimal graphs over closed starlike domains in $\l^3$ and  $\r^3,$ respectively, are relative parabolic. 
\end{abstract}

\section{Introduction}

A Riemann surface $\rb$ with non empty boundary is said to be relative parabolic if bounded harmonic functions are determined by their boundary values. This is equivalent to the existence of a positive proper superharmonic function on the surface. If $\partial (\rb)=\emptyset,$ the surface is said to be parabolic if it does not carry any non constant positive superharmonic function. See \cite{ahlfors} and \cite{grigor} for a good setting.

The conformal type problem has strongly influenced  the modern theory of minimal surfaces in the Euclidean space $\r^3$ (see for instance 
\cite{osserman}, \cite{fang}, \cite{c-k-m-r} and \cite{lop-per}, among others). 

In the Lorentzian ambient,  Calabi \cite{calabi} proved that complete maximal surfaces without singularities in $\l^3$ are spacelike planes. However, there is a vast family of complete  maximal surfaces with singularities. We emphasize the family of complete maxfaces of finite type (see \cite{um-ya}), all of them  of parabolic type by Huber's theorem \cite{hub}. 

In this paper we have obtained some parabolicity criteria for maximal and minimal surfaces {\em with non empty boundary}. Our main result asserts:
\begin{quote}
{\bf Theorem} {\em Let $X: \mb \to \l^3$ be a conformal maximal immersion. If the map $p \in \mb \mapsto \langle X(p),X(p) \rangle$ is eventually positive\footnote{That is, positive outside a compact set.}  and proper, where $\langle , \rangle$ is the Lorentzian metric, then $\mb$ is relative parabolic.

In particular, properly immersed maximal surfaces contained in the conical region $\{(x,x_3) \in \c \times \r \equiv \r^3 \;:\; |x_3| \leq \|x\| \tan (\alpha)\},$ with $\alpha \in ]0,\pi/4[,$  are relative parabolic.}
\end{quote}
The preceding theorem holds even if the immersion is singular at boundary points. In particular, it applies for maximal graphs in $\l^3$ over proper regions of $\{x_3=0\}$  with isolated singularities. In this context, relative parabolicity means that bounded harmonic functions are uniquely determined by their values at the boundary of the graph and the interior isolated singularities. To be more precise, we have obtained the following Corollary:

\begin{quote}
{\bf Corollary} {\em Let $\sb$ be a maximal surface with a closed  set of interior isolated singularities, and suppose $\sb$ is a graph over a closed starlike region in the plane $\{x_3=0\}.$ 

Then $\sb$ is relative parabolic. 

In particular, entire maximal graphs with a closed set of isolated singularities are relative parabolic.}
\end{quote}
 
A well-known standing conjecture by Meeks asserts that any minimal graph in $\r^3$ over a proper domain is  relative parabolic. Maximal and minimal surfaces are naturally interlaced via the Weierstrass representation. We exploit this connection to obtain the following parabolicity criterium for minimal graphs in $\r^3.$ 

\begin{quote}
{\bf Corollary} {\em Let $\sb$ be a minimal graph over a closed starlike region in the plane $\{x_3=0\}.$ 

Then $\sb$ is relative parabolic.}
\end{quote}

\section{Notations and Preliminary results} \label{sec:prelim}

The Euclidean metric and norm in $\r^3$  will be denoted by $\langle, \rangle_0$ and $\| \cdot \|_0,$ respectively. 

We denote by $\l^3$ the three dimensional
Lorentz-Minkowski space $(\r^3,\langle , \rangle),$ where $\langle ,
\rangle=dx_1^2+dx_2^2-dx_3^3.$ The Lorentzian {\em "norm"} is given by $\|(x_1,x_2,x_3)\|^2=x_1^2+x_2^2-x_3^2.$
We say that a vector ${\bf v} \in \r^3- \{ {\bf 0} \}$ is spacelike,
timelike or lightlike if
$\|v\|^2:=\langle {\bf v}, {\bf v} \rangle$ is positive, negative or zero,
respectively. The vector
${\bf 0}$ is spacelike by definition. When $v$ is spacelike, $\|v\|$ is chosen non negative.   A plane in $\l^3$ is spacelike,
timelike or lightlike if the induced metric is Riemannian, non degenerate
and indefinite or degenerate, respectively. We call ${\cb}_0:=\{x \in \l^3 \;:\; \|x\|=0\}$ the {\em light cone} of $\l^3$ and denote by $\mbox{Ext}(\cb_0):=\{x\in\l^3\;:\; \|x\|^2>0\}.$

Throughout this paper, ${\pi_0}:\r^3 \to \{x_3=0\}$ will denote the (Lorentzian or Euclidean) orthogonal projection.

A smooth curve $\alpha$ in $\l^3$ is said to be spacelike (resp., lightlike, timelike) if its tangent vector field is spacelike (resp., lightlike,timelike).\\

In what follows,  $\mb$ will denote a differential surface. We allow $\partial(\mb) \neq \emptyset,$ in which case we assume $\mb \subseteq \mb',$ where $\mb'$ is an open surface,  and $\partial(\mb)$ is the topological frontier $\mbox{Fr}(\mb) \subset \mb'.$
We also suppose that $\partial(\mb)$ consists of a proper (i.e., without accumulation) family of pairwise disjoint  closed curves in $\mb'$ at least piecewise ${\cb}^1,$ and some of them could be isolated points.  By definition, $\mbox{Int}(\mb)=\mb-\partial(\mb).$
A map  $X:\mb \longrightarrow \l^3$ is said to be smooth if  it  is the restriction of a smooth map on $X':\mb' \to \l^3.$ 

We say that a property holds {\em eventually} on a differentiable surface if it is valid outside a compact subset.

\begin{definition}
We say that a smooth map $X:\mb\to\l^3$ is spacelike if the tangent plane at any {\em interior} point is spacelike, that is to say, the induced metric on $\mbox{Int}(\mb)$ is Riemannian.
In this case, $\sb=X(\mb)$ is said to be a spacelike surface in $\l^3.$ 

A point $p \in \partial (\mb)$ is said to be singular if the tangent plane $T_p \mb$ with the induced metric is not Riemannian. A curve $\alpha \subseteq \partial (\mb)$ is singular if all its points are singular.

\end{definition}

We call $\h^2 = \{ (x_1,x_2,x_3) \in \r^3 \;:\;
x_1^2+x_2^2-x_3^2=-1\}$ the hyperbolic sphere in $\l^3$ of constant
intrinsic curvature $-1.$ Note that $\h^2$ has two connected  components $\h^2_+:=\h^2 \cap \{x_3 \geq 1\}$ and $\h^2_-:=\h^2 \cap \{x_3 \leq -1\}.$ The stereographic
projection $\sigma$ for $\h^2$ is defined as follows:
$$\sigma:\c \cup \{\infty\} - \{|z|=1\} \longrightarrow \h^2 \,; \; z \rightarrow
\left(\frac{-2 \mbox{Im} (z)}{1-|z|^2}, \frac{2 \mbox{Re} (z)}{1-|z|^2},
\frac{|z|^2+1}{|z|^2-1} \right),$$ where $\sigma(\infty)=(0,0,1).$\\

If $X:\mb\to\l^3$ is a spacelike immersion, the locally well defined Gauss map $N_0$ of $X$ assigns to each point of $\mbox{Int}(\mb)$ a point of $\h^2.$  A connection argument gives that $N_0$ is globally well defined and $N_0(\mbox{Int}(\mb))$ lies, up to a Lorentzian isometry, in $\h^2_-.$ This means that $\mb$ is orientable.

A Riemann surface $\rb$ is a complex manifold of dimension $1.$ As above, we allow that $\partial (\rb) \neq \emptyset,$ and in this case we always suppose that $\rb \subseteq \rb',$ where $\rb'$ is an open Riemann surface,  and $\partial(\rb) \subset \rb'$  consists of a proper family of pairwise disjoint ${\cb}^0$ closed curves in $\rb'$ (some of them could be isolated points).
By definition, a function or $1$-form on  $\rb$ is said to be harmonic (holomorphic,...) if it the restriction of a harmonic (holomorphic,...) function or $1$-form on $\rb'.$ Likewise, we define the concept of conformal map $X:\rb \to \l^3.$  

We need the following definition:
\begin{definition}[\cite{ahlfors}]
A Riemann surface $\rb$ with non empty boundary is said to be relative parabolic if the only bounded harmonic function $f$ vanishing on $\partial (\rb)$ is the constant function $f=0.$ This is equivalent to say that $\rb$ admits a proper positive  superharmonic function.

If $\partial (\rb)=\emptyset,$ $\rb$ is said to be parabolic if  positive superharmonic functions are constant.
\end{definition}

Closed regions of parabolic or relative parabolic Riemann surfaces  are relative parabolic  \cite{grigor},\cite{ahlfors}.

\subsection{Maximal surfaces}
Let $\mb$ be a Riemann surface.
A conformal map $X:\mb \longrightarrow \l^3$ is said to be a maximal immersion if $X$ is spacelike and $X|_{\mbox{Int}(\mb)}$ has null mean curvature. In this case, $\sb=X(\mb)$ is said to be a maximal surface in $\l^3.$ 

If $X:\mb \longrightarrow \l^3$ is maximal, the map $g \df \sigma^{-1} \circ N_0$ is  meromorphic on $\mbox{Int}(\mb).$ Moreover, there exists a holomorphic 1-form in $\mbox{Int}(\mb),$ such that 
\begin{equation}\label{eq:wei}
\phi_1=\frac{i}{2} (\frac{1}{g}-g)\phi_3, \quad \phi_2=-\frac{1}{2} (\frac{1}{g}+g) \phi_3
\end{equation}
are holomorphic on $\mbox{Int}(\mb),$ $\Phi=(\phi_1,\phi_2,\phi_3)$ never vanishes on $\mbox{Int}(\mb).$ We suppose that $(g,\phi_3)$ are meromorphic data on $\mb$ (i.e., they extend meromorphically beyond $\partial(\mb)$).

Up to a translation, the immersion is given by $X= \mbox{Re}  \int_{P_0} ( \phi_1,\phi_2,\phi_3),$ $P_0\in\mb.$ 

The induced Riemannian metric $ds^2$ on $\mbox{Int}(\mb)$ is given by
$ds^2=  |\phi_{1}|^2 +|\phi_{2}|^2- |\phi_{3}|^2 =\left( \frac{|\phi_3|}{2}
(\frac{1}{|g|}-|g|) \right)^2.$
Since $X$ is spacelike, then $|g| \neq 1$ on $\mbox{Int}(\mb),$ and up to a Lorentzian isometry, we always assume $|g|<1.$

We call $(\mb, \phi_1,\phi_2,\phi_3)$ (or simply
$(\mb,g,\phi_3)$)  the Weierstrass representation of $X$. For more details see, for instance, \cite{kobayashi}.

\begin{remark} \label{re:cambio}
The transformation $(\mb,\phi_1,\phi_2,\phi_3) \rightarrow (\mb,\phi_1,\phi_2,i\phi_3)$ converts Weierstrass
data of maximal surfaces in $\l^3$ into Weierstrass data of minimal surfaces in $\r^3,$ and vice versa. Moreover, the composition of the Gauss map of each surface with the corresponding stereographic projection leads to the same meromorphic map $g.$
For more details about  minimal surfaces, see \cite{osserman}.
\end{remark}

Let $X:\mb \longrightarrow \l^3$ be a conformal maximal immersion. Isolated points and loops in $\partial (\mb)$ could determine two kinds of isolated singularities in the surface $X(\mb) \subset \l^3:$ {\em branch} points and {\em lightlike} singularities. By definition, a branch point of $X(\mb)$ is the image under $X$ of an isolated singular point  $p_0 \in \partial (\mb).$ An isolated lightlike singularity in $X(\mb)$ correspond to a singular loop $\gamma\subseteq\partial(\mb)$ whose image $X(\gamma)$ is a single point. Obviously,  not all  singular loops must determine isolated singularities. 

>From the conformal  point of view, a tubular open neighborhood in $\mb$ of a singular loop $\gamma$ determining a lightlike singularity is  biholomorphic to $\{0<r<|z|\leq1\} \subset \c$ (where $\gamma \equiv \{|z|=1\}$),  while the one of a branch point $p_0$ is conformally equivalent to the open unit disc $\d$ ($p_0 \equiv 0$). Moreover, if $p_0 \in \mb$ is a branch point, the Weierstrass data extend analitically to $p_0$ and satisfy $|g(p_0)|<1$, $\phi_3(p_0)=0.$  If $\gamma$ is a singular loop determining a lightlike singularity, $\phi_3$ and $g$ also extend analytically to $\gamma,$  and in this case $|g(p)|=1,$ $p \in \gamma.$

If $p_0 \in \partial (\mb)$ is a branch point, the map $X\circ {\pi_0}:\mb\to\{x_3=0\}$ presents a topological branch point at $p_0,$ and in particular $X$ is not an embedding locally around $p_0.$ 
The geometry of isolated lightlike singularities in $X(\mb)$ is also well known. Indeed, consider the quotient surface $\hat{\mb}$ obtained by identifying all the points of the singular curve $\gamma$ with a single point $\hat{\gamma}$ and induce $\hat{X}:\hat{\mb}\to\l^3$ in the natural way. Then $\hat{X}\circ {\pi_0}:\hat{\mb}\to\{x_3=0\}$ is locally around  $\hat{\gamma}$ either an embedding  or a branched covering with branch point $\hat{\gamma}.$  In the first case, $\hat{X}$ is asymptotic to a half lightcone with vertex at $\hat{\gamma},$ and the point $\hat{X}(\hat{\gamma})=X(\gamma)$ is called a 
{\em conelike singularity}. For more details we refer to the works \cite{kobayashi}, \cite{ecker}, \cite{kly-mik} or  \cite{f-l-s}, among others.\\

If we label $\partial(\mb)_0$ as the union of brach points and singular loops in  $\partial(\mb)$ determining isolated singularities of $X(M),$ we define  $$\mbox{Int}(X(M)):=X(\mbox{Int}(M) \cup \partial(\mb)_0),$$ and observe that $\mbox{Int}(X(M))$ is a {\em branched} surface in $\r^3.$ Note that $X(\mbox{Int}(\mb))\subset \mbox{Int}(X(M)),$ but they could not coincide.

		%%%%%%%%%%%%%%%%%%%%%%%%%%%%%%%%%%%%%%%%%%%%%%%%%%%%%%%%%%%%%%%%%%%%%%%%%%%%%%%%%%%%%%%%%%%%%%%%%%%
\section{Parabolicity of maximal surfaces in $\l^3$}

We are going to prove the following theorem: 

\begin{theorem} \label{th:parabo}
Let $X:\mb \to \l^3$ be a conformal maximal immersion, where $\partial (\mb) \neq \emptyset,$ and suppose that the Lorentzian norm
$$n:\mb \to \r, \quad n(p)=\|X(p)\|^2$$
is eventually positive and proper. 

Then $\mb$ is relative parabolic.  
\end{theorem}
\begin{proof}
Consider the compact set $K=\{p\in\mb\;:\;\|X(p)\|^2\leq 2\},$ and note that $\mb$ is relative parabolic if and only if $\mb-\rm{Int}(K)$ is relative  parabolic (see \cite{ahlfors}, \cite{grigor} for details). Therefore, it suffices to check that the proper positive function $h:\mb-\rm{Int}(K)\to\r$ given by $$h(p)=\log(\|X(p)\|^2)$$ is superharmonic. Take an isothermal parameter $z=u+iv$ on $\mb.$ The corresponding conformal parameterization $X(u,v)$ satisfies  $<X_u, X_v>=0$ and $<X_u,X_u>=<X_v,X_v>=\lambda^2.$ Furthermore, from the maximality, the map $X$ is harmonic, and so:
$$\Delta h:=h_{uu}+h_{vv}=-4\left(
\frac{<X,X_u>^2+<X,X_v>^2}{<X,X>^2}-\frac{\lambda^2}{<X,X>}
\right)$$
Since  $\{X_u,X_v\}$ is an isothermal basis of a spacelike plane, we get $$X=\frac{1}{\lambda^2}\left(<X,X_u>X_u+ <X,X_v>X_v\right)-<X,N_0> N_0,$$ where $N_0$ is the normal vector. Hence, $<X,X>=\frac{1}{\lambda^2}\left(<X,X_u>^2+ <X,X_v>^2\right)-<X,N_0>^2,$ which proves that $\Delta h=-4 \lambda^2 \frac{<X,N_0>^2}{<X,X>^2}\leq 0$ and concludes the proof.
\end{proof}

\begin{corollary} \label{co:parabo}
Let $X:\mb \to \l^3$ be a proper conformal maximal immersion, where $\partial (\mb) \neq \emptyset,$ and suppose that $X(M)$ eventually lies in the conical region ${\wb}_\alpha=\{(x,x_3) \in \c \times \r \equiv \r^3 \;:\; |x_3| \leq \|x\| \tan (\alpha)\},$ where $\alpha \in ]0,\pi/4[.$  

Then $\mb$ is relative parabolic.  
\end{corollary}

\begin{proof} 
As  $X(\mb)$ is eventually contained in ${\wb}_\alpha,$  then we have
$n(p)\geq \|\pi_0(X(p))\|_0^2 (1-\mbox{tg}(\alpha)).$
Since $\alpha \in ]0,\pi/4[$ and $X$ is proper the Lorentzian norm $n:\mb \to \r$ is eventually positive and proper. The corollary follows from Theorem \ref{th:parabo}.
\end{proof}

\begin{theorem} \label{th:star}
Let $X:\mb \to \l^3,$ $\partial (\mb) \neq \emptyset,$  be a conformal proper maximal immersion\footnote{Recall that $\partial(\mb)$ could contain singular points.} and suppose that $\sb:=X(\mb)$  is a graph
over a closed starlike region $\Omega \subset \{x_3=0\} \equiv \c$  centered at $(\pi_0 \circ X)(p_0),$ where $p_0 \in \mbox{Int}(\mb).$

Then, $\mb$ is relative parabolic.
\end{theorem}
The proof of this theorem is a straightforward consequence of Theorem \ref{th:parabo} and the following Lemma.
\begin{lemma}\label{lem:posiciones}
Let $X:\mb \to \l^3$ be a proper spacelike immersion, and suppose that $\sb:=X(\mb)$  is a graph
over a closed starlike region $\Omega \subset \{x_3=0\} \equiv \c$ centered at the origin $0 \in \Omega \cap X(\mbox{Int}(\mb)).$

The following statements hold:
\begin{enumerate}[(i)]
\item ${\sb}-\{0\}$ is contained in  $\mbox{Ext}(\C_0).$
\item  If we write ${\sb}=\{(z,u(z)) \,:\; z\in \Omega\}$ and take $\theta\in[0,2\pi],$ then the function $f_\theta(t):={\rm dist}\left((t e^{i \theta},u(t e^{i \theta})),\C_0)\right)$ is positive and non decreasing in $]0,t_\theta[,$ where $t_\theta:=\mbox{Sup}\{t\in\r\;:\;te^{i\theta}\in\Omega\}\in]0,+\infty]$ (here {\em dist} means Euclidean distance).
\item The Lorentzian norm $n:{\mb} \to \r,$ $n(p)=\|X(p)\|^2$ is  non negative and proper.
\end{enumerate}
\end{lemma}
\begin{proof}
For each $\theta\in[0,2\pi],$ we label $\Pi_\theta$ as the half plane in $\l^3\equiv\c\times\r$ given by $\{(\lambda e^{i \theta},\mu) \,:\, \lambda \geq 0,\, \mu \in \r\}.$ Then,
$$\mbox{Int}({\sb})\cap\Pi_\theta=\{\rho_\theta(t):=(t e^{i \theta},u_\theta(t))\;:\;t\in[0,t_\theta[\},\quad {\rm{where}}\; u_\theta(t):=u(t e^{i\theta})$$
Since $X$ is spacelike, it is not hard to see that $u_\theta$ belongs to the Sobolev space ${\wb}^{{1},\infty}(I),$ for any compact subinterval $I \subset [0,t_\theta[.$ Furthermore, $|u'_\theta(t)| \leq 1$ on $[0,t_\theta[$ and $|u'_\theta(t)| < 1$ on regular (spacelike) points $\rho_\theta(t).$

Integrating from $t=0,$ we get that $|u_\theta(t)| \leq t$ and so ${\sb}\cap\Pi_\theta$ is contained in $\overline{\mbox{Ext}(\C_0)}.$ Since $0$ is a regular point, then $|u'_\theta(0)| < 1,$ and so,   $|u_\theta(t)| < t,$ $t > 0.$ This obviously implies that $({\sb}-\{0\})\cap\Pi_\theta \subset \mbox{Ext}(\C_0),$  $\theta \in[0,2\pi],$ 
and proves $(i).$

For $(ii),$ notice that  
$$f_\theta(t)={\rm{dist}}(\rho_\theta(t),\C_0)=\frac{1}{\sqrt{2}}\mbox{Min}\{|t-u_\theta(t)|,|t+u_\theta(t)|\}, \quad t\in[0,t_\theta[.$$
Since  $\|\rho_\theta(t)\|>0$ and $\|\rho_\theta '(t)\| \geq 0,$ $t \in ]0,t_\theta[,$ it is easy to check that $f_\theta(t)$ is positive and non decreasing in $]0,t_\theta[.$ 

To see $(iii),$ observe that the Lorentzian spheres $\h_s:=\{x\in\l^3 \;:\; \|x\|=s\},$ $s>0,$ are asymptotic to the light cone $\C_0$ in the following Euclidean sense: $\lim_{k \to\pm\infty} (r_{k,s}-r_k) = 0,$ where $r_{k,s}:=\sqrt{k^2+s^2}$ and $r_k:=|k|$ are the radii of the concentric Euclidean circles $\h_s\cap\{x_3=k\}$ and $\C_0\cap\{x_3=k\}$ respectively. 

Then it suffices  to check that
there exists $\epsilon>0$ such that $\mbox{dist}\left( X(p),\C_0\right)\geq \epsilon$ eventually in $\mb.$

Indeed, take $0<\delta<\mbox{dist}(0,\partial(\Omega))$ and observe that for $t>\delta$ and $\theta\in[0,2\pi],$
$$f_\theta(t)\geq f_\theta (\delta)\geq \epsilon:=\mbox{Min}\{ f_\theta(\delta)\;:\; \theta\in[0,2\pi]\}>0$$
which proves $(iii).$
\end{proof}

\begin{corollary}[Calabi's Theorem \cite{calabi}]\label{co:calabi}

Let $X:\mb\to\l^3$ a complete conformal maximal immersion, where $\partial(\mb)=\emptyset.$

Then $X(\mb)$ is a spacelike plane. 
The result remains valid if we substitute the hypothesis of completeness for the one of properness.
\end{corollary}
\begin{proof}
It is well known that complete (or proper) spacelike surfaces without boundary are graphs over $\{x_3=0\}.$ Without loss of generality, assume that  $0 \in X(\mb)$ and let $D \subset \{x_3=0\}$ denote an open disc centered at the origin.  By Lemma \ref{lem:posiciones}, the Lorentzian norm is non negative and proper on $\mb.$ By Theorem \ref{th:parabo}, the Riemann surface with boundary $\mb_0:=\mb-({\pi_0} \circ X)^{-1}(D)$ is 
relative parabolic, and since $\partial ( \mb_0)$ is compact, $\mb$ is conformally equivalent to $\c.$ As the the stereographic projection $g$ of the Gauss map of $X$ is a bounded holomorphic function (recall that we can suppose $|g|<1$), then we infer that $g$ is constant. The corollary follows immediately.
\end{proof}

\begin{corollary}\label{co:isolated}

Let $X:\mb\to\l^3$ a proper confomal maximal immersion, and suppose $\partial(\mb)$ consits of a countable union of  pairwise disjoint   loops determining conelike singularities.
Then $\mb$ is relative parabolic. 
\end{corollary}

\begin{proof}
Since all the singularities of $X(\mb)$ are of conelike type and the immersion $X$ is proper, $X(\mb)$ is an entire graph over any spacelike plane (see \cite{f-l}). By Theorem \ref{th:star} $\mb$ is relative parabolic.
\end{proof}

The surface $\mb$ in the previous Corollary can be biholomorphic to either $\c-\cup_{n\in\n} D_n$ or $\d-\cup_{n\in\n} D_n,$ where $\{D_n\;:\;n\in\n\}$ is a family of open discs with pairwise disjoint closures (see \cite{circulardomains}). 
However, if in addition $X(\mb)$ is invariant under a non trivial discrete group $G$ of ambient isometries acting properly and freely on $\l^3,$ and the quotient surface $\mb/G$ has a finite number of boundary components, it is possible to see that in fact $\mb\equiv\c-\cup_{n\in\n} D_n$ (see \cite{f-l}).

\subsection{Some consequences for minimal surfaces in $\r^3$}

The preceding results can be applied to obtain some parabolicity criteria for minimal surfaces.

\begin{corollary}
Let $Y\equiv(Y_1,Y_2,Y_3):\mb \to \r^3$ be a conformal minimal immersion, where $\partial (\mb) \neq \emptyset,$ and suppose that
${\pi_0} \circ Y:\mb \to \{x_3=0\}$ is proper.
Let $N \equiv (N_1,N_2,N_3):\mb \to \s^2$ denote the Gauss map of $Y,$ and assume that the real harmonic $1$-form
$\psi:=N_2 dY_1-N_1 dY_2$ is exact.

If there are $p_0 \in \mbox{Int}(\mb),$ $\epsilon >0$ and $C \in \r$ such that 
$$|\int_{p_0}^p \psi+C|\leq \|({\pi_0} \circ Y)(p)\|_0-\epsilon, \; \mbox{for all} \; p \in \mb,$$
then $\mb$ is relative parabolic.
\end{corollary}
\begin{proof}
It is easy to check that $\psi=dX_3,$ where $X_3$ denotes the harmonic conjugate of $Y_3.$ 

Up to a translation, we can suppose that $X_3(p)=\int_{p_0}^p \psi+C,$ and so  $|X_3(p)|\leq \|({\pi_0} \circ Y)(p)\|_0-\epsilon.$ Then, it is not hard to check that the maximal immersion  given by $X:=(Y_1,Y_2,X_3)$ (see Remark \ref{re:cambio}) satisfies the hypothesis of Theorem \ref{th:parabo}. This concludes the proof.
\end{proof}

\begin{corollary}
Let $Y\equiv(Y_1,Y_2,Y_3):\mb \to \r^3$ be a conformal proper minimal immersion, and suppose that $\sb:=Y(\mb)$  is a graph
over a closed starlike region $\Omega \subset \{x_3=0\} \equiv \c$  centered at  $(\pi_0 \circ Y)(p_0),$ where  $p_0 \in \mbox{Int}(\mb).$

Then, $\mb$ is relative parabolic.
\end{corollary}
\begin{proof}
Since $\mb$ is simply connected, the harmonic conjugate $X_3$ of $Y_3$ is well defined. Moreover, the maximal immersion
$X:=(Y_1,Y_2,X_3)$ satisfies the hypothesis of Theorem \ref{th:star}. The corollary follows immediately.
\end{proof}

%%%%%%%%%%%%%%%%%%%%%%%%%%%%%%%%%%%%%%%%%%%%%%%%%%%%%%%%%%%%%%%%%%%%%%%%%%%%%%%%%%%%%%%%%%%%%%%%%%

%%%%%%%%%%%%%%%%%%%%%%%%%%%%%%%%%%%%%%%%%%%%%%%%%%%%%%%%%%%%%%%%%%

{\bf ISABEL FERNANDEZ, FRANCISCO J. LOPEZ,} \newline
Departamento de Geometr\'{\i}a y Topolog\'{\i}a \newline
Facultad de Ciencias, Universidad de Granada \newline
18071 - GRANADA (SPAIN) \newline
e-mail:(first author) isafer@ugr.es, (second author) fjlopez@ugr.es

\end{document}